\documentclass[graybox]{svmult}
\usepackage{type1cm}        
\usepackage{makeidx}         
\usepackage{graphicx}        
\usepackage{multicol}        
\usepackage[bottom]{footmisc}

\usepackage{newtxtext}       %
\usepackage{newtxmath}       

\makeindex             

\newcommand\scalemath[2]{\scalebox{#1}{\mbox{\ensuremath{\displaystyle#2}}}}
\usepackage{tikz}
\usepackage{algorithm,algorithmic}

\begin{document} 

\title*{Should multilevel methods for discontinuous Galerkin
    discretizations use discontinuous interpolation operators?}
\titlerunning{Should multilevel methods for DG discretizations use
  a discontinuous interpolation?}

\author{Martin J. Gander and José Pablo Lucero Lorca}

\institute{Martin J. Gander \at University of Geneva,
  \email{martin.gander@unige.ch} \and José Pablo Lucero Lorca \at
  University of Colorado at Boulder, \email{pablo.lucero@colorado.edu}}

\maketitle

\abstract*{Multi-level preconditioners for Discontinuous Galerkin (DG)
discretizations are widely used to solve elliptic equations, and a
main ingredient of such solvers is the interpolation operator to
transfer information from the coarse to the fine grid. Classical
interpolation operators give continuous interpolated values, but since
DG solutions are naturally discontinuous, one might wonder if one
should not use discontinuous interpolation operators for DG
discretizations. We consider a discontinuous interpolation operator
with a parameter that controls the discontinuity, and determine the
optimal choice for the discontinuity, leading to the fastest solver
for a specific 1D symmetric interior penalty DG discretization model
problem. We show in addition that our optimization delivers a
perfectly clustered spectrum with a high geometric multiplicity, which
is very advantageous for a Krylov solver using the method as its
preconditioner. Finally, we show the applicability of the optimal
choice to higher dimensions.  }

\section{Discontinuous Interpolation for a Model Problem}

Interpolation operators are very important for the construction of a
multigrid method.  Since multigrid's inception by Fedorenko
\cite{Fedorenko1964}, interpolation was identified as key, deserving
an entire appendix in Brandt's seminal work \cite{Brandt1977}:
\emph{'[...] even a small and smooth residual function may produce
large high-frequency residuals, and significant amount of
computational work will be required to smooth them out.'}

With the advent of discontinuous Galerkin (DG) discretizations
\cite{Arnold1982}, the problem of choosing an interpolation becomes
particularly interesting. A good interpolation operator will not
produce undesirable high frequency components in the residual. For DG
discretizations however, high frequencies seem to be inherent. In a
wider context, the choice of restriction and prolongation operators
defines the coarse space itself when using an inherited (Galerkin)
coarse operator, and convergence of multigrid algorithms with
classical restriction and interpolation operators for DG
discretizations of elliptic problems cannot be independent of the
number of levels if inherited bilinear forms are considered
\cite{AntoniettiSartiVerani2015}. In 1D, the reason for this was
identified in \cite[\S4.3]{GanderLucero2020}): the DG penalization is
doubled at each coarsening, causing the coarse problem to become
successively stiffer.

A simple classical interpolation operator is linear interpolation: in
1D one takes the average from two adjacent points in the coarser grid
and sets the two DG degrees of freedom at the midpoint belonging to
the fine mesh to this same value, therefore imposing continuity at
that point.  But why should continuity be imposed on the DG
interpolated solution on the fine mesh?  Is it possible to improve the
solver performance with a discontinuous interpolation operator?

Convergence of two-level methods for DG discretizations has been
analyzed for continuous interpolation operators using classical
analysis, see \cite{FengKarakashian2001} and references therein, and
also Fourier analysis \cite{Hemker1980, Hemker2003, Hemker2004,
GanderLucero2020}.  We use Fourier analysis here to investigate the
influence of a discontinuous interpolation operator on the two level
solver performance.  We consider a symmetric interior penalty
discontinuous Galerkin (SIPG) finite element discretization of the
Poisson equation as in \cite{Arnold1982} on a 1D mesh as shown in
Fig. \ref{fig:Mesh}.
\begin{figure}[t]
  \centering
  \begin{tikzpicture}[scale=1.]
    \draw[-] (0.0,0.0) -- (10.5, 0.0);
    \draw[-] (0.0,0.1) -- (0.0,-0.1) node[below]{$x_{1}^+ = 0$} ;
    \draw[-] (1.5,0.1) -- (1.5,-0.1) node[below]{$x_{1}^- x_{2}^+$};
    \draw[-] (3.0,0.1) -- (3.0,-0.1) node[below]{$\dots$};
    \draw[-] (4.5,0.1) -- (4.5,-0.1) node[below]{$x_{j-1}^-x_{j}^+$};
    \draw[-] (6.0,0.1) -- (6.0,-0.1) node[below]{$x_{j}^- x_{j+1}^+$};
    \draw[-] (7.5,0.1) -- (7.5,-0.1) node[below]{$\dots$};
    \draw[-] (9.0,0.1) -- (9.0,-0.1) node[below]{$x_{J-1}^- x_{J}^+$};
    \draw[-] (10.5,0.1) -- (10.5,-0.1) node[below]{$x_{J}^- = 1$};
  \end{tikzpicture}
  \caption{Mesh for our DG discretization of the Poisson
    equation.}\label{fig:Mesh}
\end{figure}
The resulting linear system reads (for details see
\cite{GanderLucero2020}) \renewcommand{\arraystretch}{1.7}
\begin{align} \label{eqn:DiscProb}
  A \boldsymbol{u} = 
  \frac1{h^2} \left(
  \begin{matrix}
    \ddots   & \ddots             & -\frac12                      &                              &                   &           \\[-0.7em]
    \ddots   & \delta_0           & 1-\delta_0                    &                   -\frac12   &                   &           \\[-0.7em]
    -\frac12 & 1-\delta_0         & \delta_0                      &                              & - \frac12         &           \\[-0.7em]
             & - \frac12          &                               & \delta_0                     &  1-\delta_0       & - \frac12 \\[-0.7em]
             &                    &                      -\frac12 & 1-\delta_0                   &          \delta_0 &    \ddots \\[-0.7em] 
             &                    &                               &                     -\frac12 &            \ddots &    \ddots
  \end{matrix}\right) \left(\begin{matrix} 
    \vdots\\[-0.7em] 
    u_{j-1}^+\\[-0.7em] 
    u_{j}^- \\[-0.7em] 
    u_{j}^+ \\[-0.7em] 
    u_{j+1}^- \\[-0.7em] 
    \vdots
  \end{matrix}\right) = \left(\begin{matrix} 
    \vdots \\[-0.7em] 
    f_{j-1}^+\\[-0.7em] 
    f_{j}^- \\[-0.7em] 
    f_{j}^+ \\[-0.7em] 
    f_{j+1}^- \\[-0.7em] 
    \vdots
  \end{matrix}\right) \eqqcolon \boldsymbol{f},
\end{align}
where the top and bottom blocks will be determined by the boundary
conditions, $h$ is the mesh size, $\delta_0 \in \mathbb{R}$ is the DG
penalization parameter, $\boldsymbol{f} = (\dots, f_{j-1}^+, f_{j}^-,
f_{j}^+, f_{j+1}^-, \dots) \in \mathbb{R}^{2J}$ is the source vector,
analogous to the solution $\boldsymbol{u}$.

The two-level preconditioner $M^{-1}$ we study consists of a
cell-wise nonoverlapping Schwarz (a \emph{cell} block-Jacobi)
smoother $D_c^{-1}$ (see
\cite{FengKarakashian2001,Dryja2016})\footnote{In 1D this is simply a
Jacobi smoother, which is not the case in higher dimensions.}, and a
new discontinuous interpolation operator $P$ with discontinuity
  parameter $c$, i.e.  \newcommand{\edots}{\rotatebox{0}{$\ddots$}}
\begin{align}\label{eqn:cellBJ}
  D_c^{-1} \boldsymbol{u} \coloneqq h^2 
  \left(\begin{array}{cccc}
          \ddots &     &     &       \\
                 & \delta_0 &     &       \\
                 &     & \delta_0 &       \\
                 &     &     & \ddots 
        \end{array}\right)^{-1}
                               \left(\begin{array}{c}
                                       \vdots  \\
                                       u_{j}^+ \\
                                       u_{j}^- \\
                                       \vdots
  \end{array}\right),\quad
  P \coloneqq \left(
  \begin{array}{cccccccc}
    1    &         &         &         \\[-1.1em]
    c    &   1-c   &         &         \\[-1.1em]
   1-c   &    c    &         &         \\[-1.1em]
         &    1    &         &         \\[-1.1em]
         &         &  \edots &         \\[-1.1em]
         &         &  \edots & \edots  \\[-1.1em]
         &         &  \edots & \edots  \\[-1.1em]
         &         &         & \edots        
  \end{array}
  \right),
\end{align}
where $c=\frac12$ gives continuous interpolation. The restriction
operator is $R \coloneqq \frac12 P^\intercal$, and we use $A_0 := R A
P$. The action of our two-level preconditioner $M^{-1}$, with one
presmoothing step and a relaxation parameter $\alpha$, acting on a
residual $g$, is given by
\begin{enumerate}
  \item compute $\boldsymbol{x}:= \alpha D_c^{-1} \boldsymbol{g}$,
\item compute $\boldsymbol{y}:= \boldsymbol{x} + P A_0^{-1} R
(\boldsymbol{g} - A \boldsymbol{x})$,
\item obtain $M^{-1}\boldsymbol{g} = \boldsymbol{y}$.
\end{enumerate}

\section{Study of optimal parameters by Local Fourier Analysis}

In \cite{GanderLucero2020} we described in detail, for classical
interpolation, how Local Fourier Analysis (LFA) can be used to block
diagonalize all the matrices involved in the definition of $M^{-1}$ by
using unitary transformations. The same approach still works with our
new discontinuous interpolation operator, and we thus use the same
definitions and notation for the block-diagonalization matrices $Q$,
$Q_l$, $Q_r$, $Q_0$, ${Q_l}_0$ and ${Q_r}_0$ from
\cite{GanderLucero2020}, working directly with matrices instead of
stencils in order to make the important LFA more accessible to our
linear algebra community. We extract a submatrix $\widetilde{A}$
containing the degrees of freedom of two adjacent cells from the SIPG
operator defined in \eqref{eqn:DiscProb},
\begin{align*}
  \widetilde{A} = \frac1{h^2}
  \begin{pmatrix}
    -\frac12 &   1-\delta_0   &   \delta_0    &     0    & -\frac12 &          &          &          \\
             & -\frac12  &     0    &    \delta_0   &  1-\delta_0   & -\frac12 &          &          \\
             &           & -\frac12 &   1-\delta_0  &   \delta_0    &     0    & -\frac12 &          \\
             &           &          & -\frac12 &    0     &    \delta_0   &   1-\delta_0  & -\frac12   
  \end{pmatrix},
\end{align*}
which we can block-diagonalize, $\widehat{A} = Q_l \widetilde{A} Q_r$,
to obtain
\begin{align*}
  \widehat{A} = \frac1{h^2}
  \scalemath{0.8}{\begin{pmatrix} 
      \delta_0+\cos\left(2 \pi (k-J/2) h \right) &                                1-\delta_0  &                                 &                                \\
      1-\delta_0                                 & \delta_0+\cos\left(2 \pi (k-J/2) h \right) &                                 &                                \\
      &                                     &       \delta_0-\cos\left(2 \pi k h \right) & 1-\delta_0                          \\
      &                                     &                                 1-\delta_0 & \delta_0-\cos\left(2 \pi k h \right)
    \end{pmatrix}}.
\end{align*}
The same mechanism can be applied to the smoother, $\widehat{D}_c =
Q_l \widetilde{D}_c Q_r = \frac{\delta_0}{h^2} I$, where $I$ is the $4
\times 4$ identity matrix, and also to the restriction and
prolongation operators,
$\widehat{R}=\frac12 {Q_l}_0 \widetilde{R} Q_r$ with
\begin{align*}
    \widehat{R}= \frac1{\sqrt{2}}
    \scalemath{0.8}{
      \begin{pmatrix}
        1+(c-1) e^{\frac{2 i \pi k}{J}} & - c e^{\frac{2 i \pi k}{J}}  & (-1)^j\left(1-(c-1) e^{\frac{2 i \pi k}{J}}\right) & (-1)^j c e^{\frac{2 i \pi k}{J}} \\
        (-1)^j c e^{-\frac{2 i \pi k}{J}} &(-1)^j\left(1+(c-1) e^{-\frac{2 i \pi k}{J}}\right) &  c e^{-\frac{2 i \pi k}{J}} &  1-(c-1) e^{-\frac{2 i \pi k}{J}} 
      \end{pmatrix}},
\end{align*}
and $P = 2 R^\intercal$, $\widehat{P} = Q_l \widetilde{P} {Q_r}_0 = 2
\widehat{R}^*$. Finally, for the coarse operator, we obtain
  $Q_0^* A_0 Q_0 = Q_0^* R A P Q_0 = Q_0^* R Q Q^* A Q Q^* P Q_0$, and
  thus $\widehat{A}_0= \widehat{R} \widehat{A} \widehat{P}$ with
\begin{align*}\label{eqn:coarseoperator}
    \widehat{A}_0\!=\! \frac1{H^2} 
    \scalemath{0.6}{\begin{pmatrix}
        \frac{1}{2} \left(c \left(4 (c-1) \delta_0-2 c+3\right)+(c-1) \cos \left(\frac{4 \pi  k}{J}\right)+2 \delta_0-1\right) & \frac{1}{2} (-1)^j \left(-(2 c-1) \left(c \left(2 \delta_0-1\right)-\delta_0+1\right) e^{\frac{4 i \pi  k}{J}}-c-\delta_0+1\right)
        \\
        \frac{1}{2} (-1)^j \left(-(2 c-1) \left(c \left(2 \delta_0-1\right)-\delta_0+1\right) e^{-\frac{4 i \pi  k}{J}}-c-\delta_0+1\right) & \frac{1}{2} \left(c \left(4 (c-1) \delta_0-2 c+3\right)+(c-1) \cos \left(\frac{4 \pi  k}{J}\right)+2 \delta_0-1\right)
      \end{pmatrix}},
\end{align*}
where $H=2h$. We notice that the coarse operator is different for $j$
even and $j$ odd; however, the matrices obtained for both cases are
similar, with similarity matrix $(-1)^j I$ where $I$ is the identity
matrix, and therefore have the same spectrum. In what follows we
assume $j$ to be even, without loss of generality.
This means that we will be studying a node that is present in both the
coarse and fine meshes.

The error reduction capabilities of our two level preconditioner
  $M^{-1}$ are given by the spectrum of the stationary iteration
operator
\begin{align*}
  E = (I - P A_0^{-1} R A)(I - \alpha D_{c}^{-1} A),
\end{align*}
and as in \cite{GanderLucero2020}, the 4-by-4 block
Fourier-transformed operator from LFA,
\begin{align*}
  \widehat{E}(k) = (I - \widehat{P}(k) \widehat{A}_0^{-1}(k)
  \widehat{R}(k) \widehat{A}(k))(I - \alpha \widehat{D}_{c}^{-1}(k)
  \widehat{A}(k)),
\end{align*}
has the same spectrum. Thus, we focus on studying the spectral radius
$\rho(\widehat{E}(k))$ in order to find the optimal choices for the
relaxation parameter $\alpha$, the penalty parameter $\delta_0$ and
the discontinuity parameter $c$.  The non zero eigenvalues of
$\widehat{E}(k)$ are of the form $\lambda_\pm := c_1 \pm
\sqrt{\frac{c_2}{c_3}}$, with
\begin{align*}
  c_1 &= \left\{\scalemath{0.8}{
  \begin{aligned}
    &\bigg\{-\alpha  \left(3 c^2 \delta _0 \left(4 \delta _0-3\right)+c \left(-12 \delta _0^2+9 \delta _0+1\right)+4 \delta _0^2-2 \delta_0-1\right)\\
    &+\delta_0 \left(c^2 \left(8 \delta_0^2-4 \delta_0-1\right)+c \left(-8 \delta_0^2+4 \delta _0+2\right)+2 \delta_0^2-1\right)\\
    &\left.+(1-c) \left(\alpha +\alpha  c \left(\delta _0-2\right)+(c-1) \delta_0\right) \cos\left(\frac{4 \pi  k}{J}\right)\right\} \bigg/ \\
    &\left(2 \delta_0^2 - 1 + \delta_0 c^2 \left(8 \delta _0^2-4 \delta _0-1\right)+ \delta_0 c \left(-8 \delta _0^2+4 \delta _0+2\right)- \delta_0 (c-1)^2 \cos\left(\frac{4 \pi  k}{J}\right)\right),
  \end{aligned}
  }\right.\\
  c_2 &= \left\{\scalemath{0.8}{
  \begin{aligned}
    &2 \alpha ^2 \bigg(16 (c-1)^2 c^2 \delta _0^4-2 (c-1)^2 \left(4 c^2+c+2\right) \delta _0-8 (c-1) c (3 (c-1) c-1) \delta_0^3\\
    &+\left(c (17 c+8) (c-1)^2+2\right) \delta _0^2+2 (c-1)^2 ((c-1) c+1)\bigg) \\
    &+4 \alpha ^2 \left(4 (c-1) c \delta _0^2-3 (c-1) c \delta _0+c+\delta _0-1\right) \left(c \left(3 (c-1) \delta _0-2 c+3\right)+\delta _0-1\right) \cos \left(\frac{4 \pi  k}{J}\right) \\
    &+2 \alpha ^2 (c-1)^2 c \left(c \left(\left(\delta _0-4\right) \delta _0+2\right)+2 \left(\delta _0-1\right)\right) \cos ^2\left(\frac{4 \pi  k}{J}\right),
  \end{aligned}
  }\right.\\
  c_3 &= \left\{\scalemath{0.8}{
  \begin{aligned}
    &\delta _0^2 \left(4 c (c-1) \delta _0-2 (1-2 c)^2 \delta _0^2+(c-1)^2\right)^2 \\
    &+2 \delta_0^2 \left(-2 \left(2 c^2-3 c+1\right)^2 \delta _0^2+4 c (c-1)^3 \delta _0+(c-1)^4\right) \cos \left(\frac{4 \pi  k}{J}\right)\\
    &+(c-1)^4 \delta_0^2 \cos^2\left(\frac{4 \pi k}{J}\right).
  \end{aligned}
  }\right.
\end{align*}
A first approach to optimize would be to minimize the spectral radius
for all frequency parameters $k$, but if we can find a combination of
the parameters $(\alpha, \delta_0, c)$ such that the eigenvalues of
the error operator do not depend on the frequency parameter $k$, then
the spectrum of the iteration operator, and therefore the
preconditioned system becomes perfectly clustered , i.e. only a few
eigenvalues repeat many times, regardless of the size of the
problem. The solver then becomes mesh independent, and the
preconditioner very attractive for a Krylov method that will converge
in a finite number of steps.

For these equations not to depend on $k$, they must be independent of
$\cos\left(\frac{4 \pi k}{J}\right)$, and to achieve this, we impose
three conditions on the coefficients accompanying the cosine, and we
deduce a combination of the parameters $(\alpha, \delta_0, c)$ which
we verify \emph{a posteriori} fall into the allowed range of values
for each parameter. Our conditions are:
\begin{enumerate}
\item Set the coefficient accompanying the cosine in the numerator
  of $c_1$ to zero.
\item Since the denominator of $c_1$ also contains the cosine, set the
  rest of the numerator of $c_1$ to zero in order to get rid of $c_1$
  entirely. Note that this requirement immediately implies an
  equioscillating spectrum, which often is characterizing the solution
  minimizing the spectral radius, see e.g. \cite{GanderLucero2020}.
\item $c_2$ and $c_3$ are second order polynomials in the cosine
  variable, if we want the quotient to be non zero and independent
    of the cosine, we need for the polynomials to simplify and for that,
    they must differ only by a multiplying factor independent of the
    cosine.  We then equate the quotient of the quadratic terms with the
    quotient of the linear terms and verify \emph{a posteriori} that
    $c_2/c_3$ becomes indeed independent of the cosine.
\end{enumerate}
These three conditions lead to the nonlinear system of equations
\begin{align*}
  \left\{
  \begin{aligned}
    &\begin{aligned}
      \alpha +\alpha  c \left(\delta _0-2\right)+(c-1) \delta_0 =& 0,
    \end{aligned}\\\quad\\
    &\begin{aligned}
      &\alpha  \left(3 c^2 \delta _0 \left(4 \delta _0-3\right)+c \left(-12 \delta _0^2+9 \delta _0+1\right)+4 \delta _0^2-2 \delta_0-1\right) = \\
      &\delta_0 \left(c^2 \left(8 \delta_0^2-4 \delta_0-1\right)+c \left(-8 \delta_0^2+4 \delta _0+2\right)+2 \delta_0^2-1\right),
    \end{aligned}\\\quad\\
    &\begin{aligned}
      &\frac{2 \alpha ^2 (c-1)^2 c \left(c \left(\left(\delta _0-4\right) \delta _0+2\right)+2 \left(\delta _0-1\right)\right)}{(c-1)^4 \delta_0^2} = \\
      &\frac{4 \alpha ^2 \left(4 (c-1) c \delta _0^2-3 (c-1) c \delta _0+c+\delta _0-1\right) \left(c \left(3 (c-1) \delta _0-2 c+3\right)+\delta _0-1\right)}
      {2 \delta_0^2 \left(-2 \left(2 c^2-3 c+1\right)^2 \delta _0^2+4 c (c-1)^3 \delta _0+(c-1)^4\right)}.
    \end{aligned}
  \end{aligned}\right.
\end{align*}
This system of equations can be solved either numerically or
symbolically. After a significant effort, the following values solve
our nonlinear system:
\begin{align*}
  c =& \text{Root of } 3 - 8 \tilde{c} + 8 \tilde{c}^2 - 8 \tilde{c}^3 + 4 \tilde{c}^4 \text{ such that } \tilde{c} \in \mathbb{R} \text{ and } 0 < \tilde{c} < 1, \\
  \delta_0 =& \text{Root of } -1 - 4 \tilde{\delta_0} + 24 \tilde{\delta_0}^2 - 32 \tilde{\delta_0}^3 + 12 \tilde{\delta_0}^4 \text{ such that } \tilde{\delta_0} \in \mathbb{R} \text{ and } 1 < \tilde{\delta_0}, \text{ and} \\
  \alpha =& \text{Root of } -1 - 40 \tilde{\alpha} + 214 \tilde{\alpha}^2 - 352 \tilde{\alpha}^3 + 183 \tilde{\alpha}^4 \text{ such that } \tilde{\alpha} \in \mathbb{R} \text{ and } 0 < \tilde{\alpha} < 1.
\end{align*}
The corresponding numerical values are approximately
$$
  c \approx 0.564604, \quad \delta_0 \approx 1.516980, \quad  \alpha \approx
  0.908154,
$$
and we see that indeed the interpolation should be discontinuous!  We
have found a combination of parameters that perfectly clusters the
eigenvalues of the iteration operator of our two level method, and
therefore also the spectrum of the preconditioned operator. Such
clustering is not very often possible in preconditioners, a few
exceptions are the HSS preconditioner in \cite{benzi2003optimization},
and some block preconditioners, see
e.g. \cite{silvester1994fast}. Furthermore, the spectrum is
equioscillating, which often characterizes the solution minimizing the
spectral radius of the iteration operator.

\section{Numerical Results}
\begin{figure}[t]
  \mbox{\includegraphics[width=0.4\textwidth]{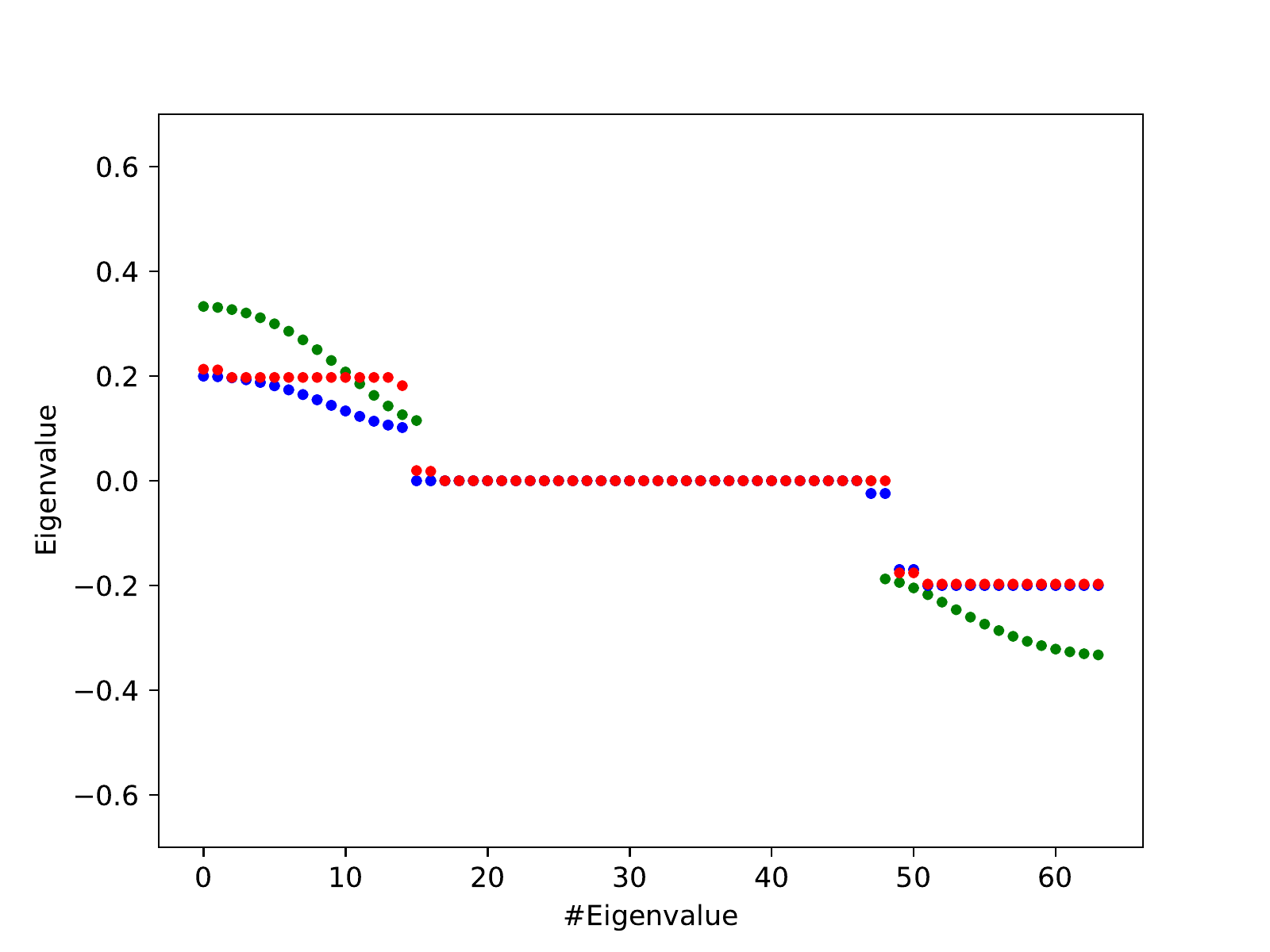}\hspace{-1em}
    \includegraphics[width=0.65\textwidth]{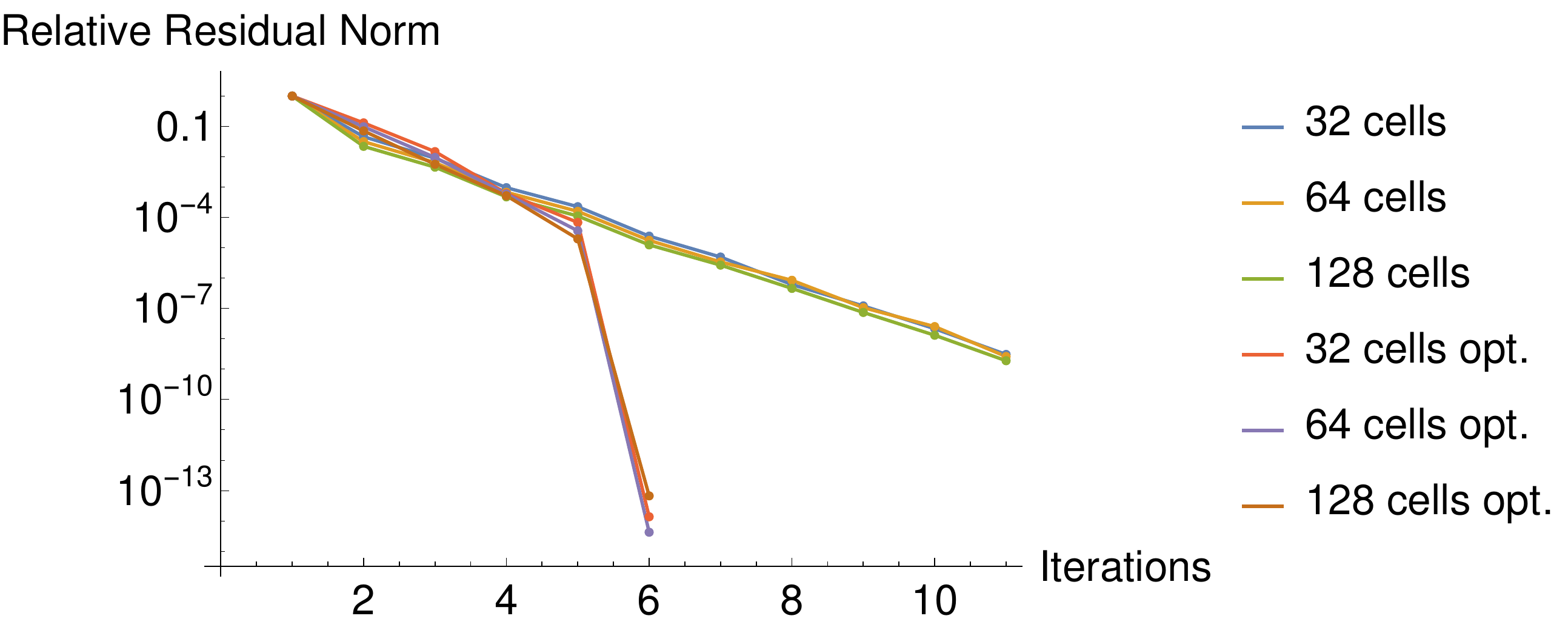}}
  \caption{Solving $-\Delta u = 1$ in 1D with Dirichlet boundary
    conditions. Left: eigenvalues of the error operator $E$, for a 32-cell
    mesh.  Green: optimizing $\alpha$ for $\delta_0=2$ (classical choice).
    Blue: optimizing $\alpha$ and $\delta_0$.  Red: optimizing $\alpha$,
    $\delta_0$ and $c$. Right: GMRES iterations for classical
    interpolation $c=0.5$, with $\delta_0=2$ and $\alpha=8/9$, and for the
    optimized clustering choice, leading to finite step convergence.}
  \label{fig:Eigenvalues}
\end{figure}

We show in Fig. \ref{fig:Eigenvalues} on the left
the eigenvalues of the iteration operator for a 32-cell mesh in 1D
with Dirichlet boundary conditions, for continuous interpolation and
$\delta_0=2$ optimizing only $\alpha$, optimizing both $\alpha$ and
$\delta_0$, and the optimized clustering choice. We clearly see the
clustering of the eigenvalues, including some extra clusters due to
the Dirichlet boundary conditions. We also note that the spectrum is
nearly equioscillating due to condition (1) and (2), which delivers
visibly an optimal choice in the sense of minimizing the spectral
radius of the error operator. With periodic boundary conditions,
the spectral radius for the optimal choice of $\alpha$, $\delta_0$ and
$c$ is $0.19732$, while only optimizing $\alpha$ and $\delta_0$ it is
$0.2$. The eigenvalues due to the Dirichlet boundary conditions are
slightly larger than $0.2$, but tests with periodic boundary
conditions confirm that then these larger eigenvalues are not present.
Refining the mesh conserves the shape of the spectrum shown in
Fig. \ref{fig:Eigenvalues} on the left, but with more eigenvalues in
each cluster, except for the clusters related to the Dirichlet
boundary conditions. Note also that since the error operator is
equioscillating around zero, the spectrum of the preconditioned system
is equioscillating around one, and since the spectral radius is less
than one, the preconditioned system has a positive spectrum and is
thus invertible.

In Fig. \ref{fig:Eigenvalues} on the right we show the GMRES
iterations needed to reduce the residuals by $10^{-8}$ for different
parameter choices and the clustering choice, for different mesh
refinements. We observe that the GMRES solver becomes exact after six
iterations for the clustering choice.

\begin{figure}[t]
  \centering
  \includegraphics[width=0.6\textwidth]{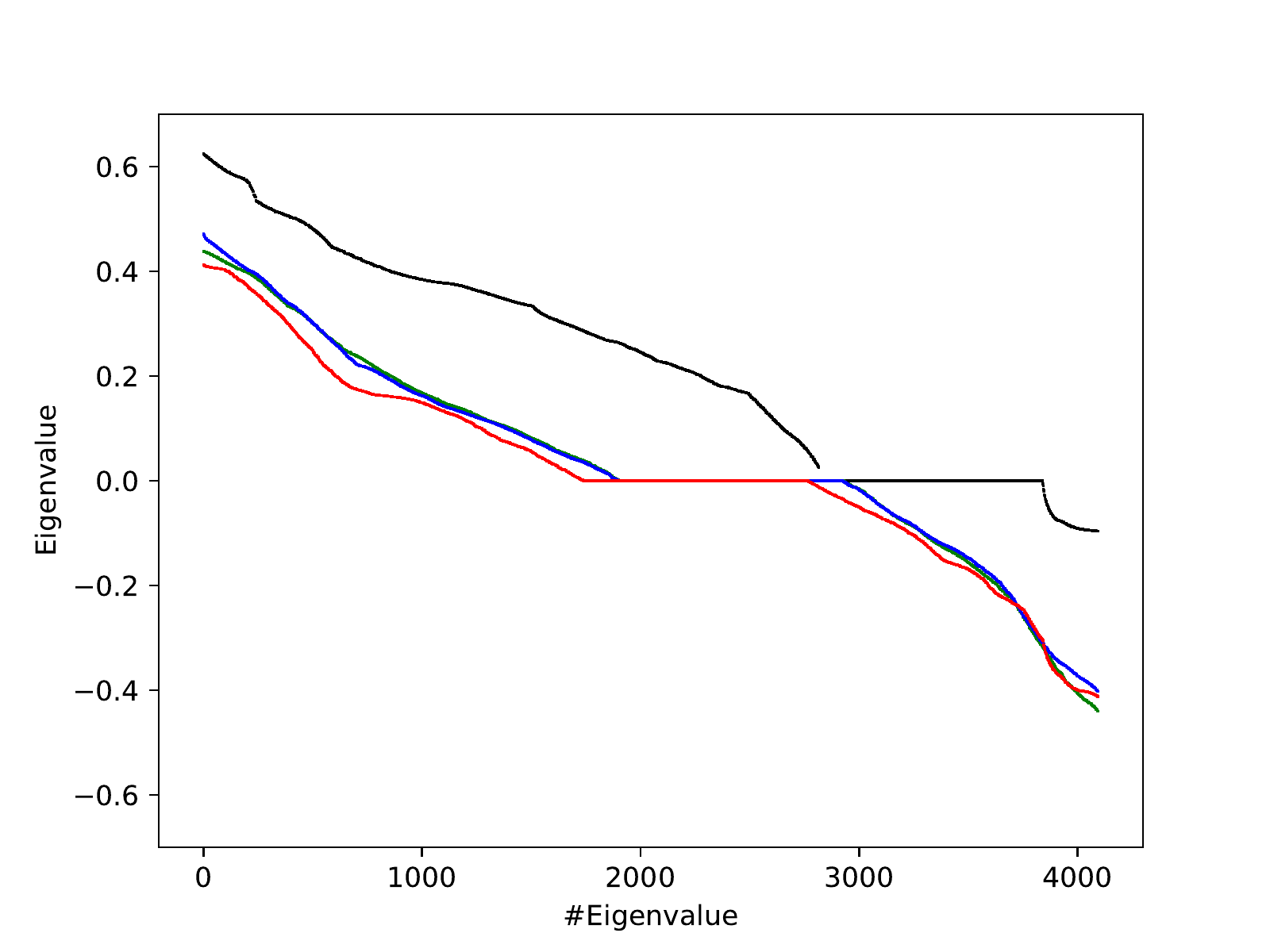}
  \caption{Spectrum of the iteration operator for a 32-by-32
    square 2D mesh.  Black: optimizing $\alpha$ for $\delta_0=2$
    (classical choice) in 1D.  Green: optimizing $\alpha$ and
    $\delta_0$ in 1D.  Blue: optimizing $\alpha$, $\delta_0$ and $c$
    in 1D.  Red: numerically optimizing $\alpha$, $\delta_0$ and $c$
    in 2D.}
  \label{fig:2Dspec}
\end{figure}

We next perform tests in two dimensions using an interpolation
operator with a stencil that is simply a tensor product of the 1D
stencil $\left(\begin{smallmatrix} 1 & 0 \\ c & 1-c \\ 1-c & c \\ 0 &
  1 \end{smallmatrix}\right) \otimes \left(\begin{smallmatrix} 1 & 0
  \\ c & 1-c \\ 1-c & c \\ 0 & 1 \end{smallmatrix}\right)$, where
$\otimes$ stands for the Kronecker product. This is very common in DG
methods where even the cell block-Jacobi matrix can be expressed as a
Kronecker sum for fast inversion.  We show in Fig.
\ref{fig:2Dspec}
the spectrum for different optimizations in two dimensions.  We
observe that the clustering is not present, however as shown in detail
in \cite{GanderLucero2020} for classical interpolation, the optimal
choice from the 1D analysis is also here very close to the numerically
calculated optimum in 2D.

\section{Conclusion}

We studied the question whether in two level methods for DG
discretizations one should use a discontinuous interpolation
operator. Our detailed analysis of a one dimensional model problem
showed that the optimization of the entire two level process indeed
leads to a discontinuous interpolation operator, and the performance
of this new two level method is superior to the case with the
continuous interpolation operator. More importantly however, the
discontinuous interpolation operator allowed us to cluster the
spectrum for our model problem, which then permits a Krylov method
with this preconditioner to become a direct solver, converging in the
number of iterations corresponding to the number of clusters in exact
arithmetic. Our numerical experiments showed that this is indeed the
case, but that when using the one dimensional optimized parameters in
higher spatial dimensions, the spectrum is not clustered any
more. Nevertheless, the 1D parameters lead to a two level performance
very close to the numerically optimized one in 2D.  We are currently
investigating if there is a choice of interpolation operator in 2D
which still would allow us to cluster the spectrum, and what the
influence of this discontinuous interpolation operator is on the
Galerkin coarse operator obtained.

\bibliographystyle{plain}
\bibliography{proc}{}

\end{document}